\documentclass[a4paper,twoside,10pt]{article}
\usepackage{amsmath}
\usepackage{amssymb}
\usepackage{amsthm}
\usepackage{amscd}
\begin{document}
%%%%%%%%%%%%%%%%%%%% Text italic %%%%%%%%%%%%%%%%%%%%%%%%%%%%
\theoremstyle{plain}
\newtheorem{thm}{Theorem}[section]
\newtheorem{theorem}[thm]{Theorem}
\newtheorem{lemma}[thm]{Lemma}
\newtheorem{corollary}[thm]{Corollary}
\newtheorem{proposition}[thm]{Proposition}
\newtheorem{addendum}[thm]{Addendum}
\newtheorem{variant}[thm]{Variant}
%%%%%%%%%%%%%%%%%%%% Text roman %%%%%%%%%%%%%%%%%%%%%%%%%%%%%
\theoremstyle{definition}
\newtheorem{notations}[thm]{Notations}
\newtheorem{question}[thm]{Question}
\newtheorem{problem}[thm]{Problem}
\newtheorem{remark}[thm]{Remark}
\newtheorem{remarks}[thm]{Remarks}
\newtheorem{definition}[thm]{Definition}
\newtheorem{claim}[thm]{Claim}
\newtheorem{assumption}[thm]{Assumption}
\newtheorem{assumptions}[thm]{Assumptions}
\newtheorem{properties}[thm]{Properties}
\newtheorem{example}[thm]{Example}

\numberwithin{equation}{section}

\newcommand{\sA}{{\mathcal A}}
\newcommand{\sB}{{\mathcal B}}
\newcommand{\sC}{{\mathcal C}}
\newcommand{\sD}{{\mathcal D}}
\newcommand{\sE}{{\mathcal E}}
\newcommand{\sF}{{\mathcal F}}
\newcommand{\sG}{{\mathcal G}}
\newcommand{\sH}{{\mathcal H}}
\newcommand{\sI}{{\mathcal I}}
\newcommand{\sJ}{{\mathcal J}}
\newcommand{\sK}{{\mathcal K}}
\newcommand{\sL}{{\mathcal L}}
\newcommand{\sM}{{\mathcal M}}
\newcommand{\sN}{{\mathcal N}}
\newcommand{\sO}{{\mathcal O}}
\newcommand{\sP}{{\mathcal P}}
\newcommand{\sQ}{{\mathcal Q}}
\newcommand{\sR}{{\mathcal R}}
\newcommand{\sS}{{\mathcal S}}
\newcommand{\sT}{{\mathcal T}}
\newcommand{\sU}{{\mathcal U}}
\newcommand{\sV}{{\mathcal V}}
\newcommand{\sW}{{\mathcal W}}
\newcommand{\sX}{{\mathcal X}}
\newcommand{\sY}{{\mathcal Y}}
\newcommand{\sZ}{{\mathcal Z}}
% Sonderbuchstaben mit Doppellinie
%\renewcommand{\chaptermark}[1]{\markboth{Ch.\,\thechapter. \, #1}{}}
%\renewcommand{\sectionmark}[1]{\markright{\thesection. \, #1}}
\newcommand{\bbA}{{\mathbb A}}
\newcommand{\bbB}{{\mathbb B}}
\newcommand{\bbC}{{\mathbb C}}
\newcommand{\bbD}{{\mathbb D}}
\newcommand{\bbE}{{\mathbb E}}
\newcommand{\bbF}{{\mathbb F}}
\newcommand{\bbG}{{\mathbb G}}
\newcommand{\bbH}{{\mathbb H}}
\newcommand{\bbI}{{\mathbb I}}
\newcommand{\bbJ}{{\mathbb J}}
\newcommand{\bbM}{{\mathbb M}}
\newcommand{\bbN}{{\mathbb N}}
\newcommand{\bbQ}{{\mathbb Q}}
\newcommand{\bbR}{{\mathbb R}}
\newcommand{\bbT}{{\mathbb T}}
\newcommand{\bbU}{{\mathbb U}}
\newcommand{\bbV}{{\mathbb V}}
\newcommand{\bbW}{{\mathbb W}}
\newcommand{\bbX}{{\mathbb X}}
\newcommand{\bbY}{{\mathbb Y}}
\newcommand{\bbZ}{{\mathbb Z}}
\newcommand{\id}{{\rm id}}
\newcommand{\rank}{{\rm rank}}

%%%%%%%%%%%%%%%%%%%%%%%%%%%%%%%%%%%%%%%%%%%%%%%%%%%%%%%%%%%%%%
\title{\textbf{On complex surfaces with $5$ or $6$ semistable singular fibers
over $\mathbb P^1$}
\thanks{Keywords: singular fiber, fibration, surface, canonical class
inequality. \  The first author would like to thank the support of
``DFG-Schwerpunktprogramm Globale Methoden in der Komplexen
Geometrie'' and the DFG-NSFC Chinese-German project ``Komplexe
Geometrie''. He is also supported by the 973 Project Foundation
and the Foundation of EMC for Doctoral Program. The third author
is partially supported by Conacyt Grant 41459-F.}}

\author{{\bf Sheng-Li Tan} \\[6pt]
Department of Mathematics, East China Normal University, \\
Shanghai
200062, P. R. of China \\
sltan@math.ecnu.edu.cn \\[8pt]
%{\it Current address}: Fb 6, Mathematik, Universit\"at Essen,
%45117 Essen, Germany \\
%sheng-li.tan@uni-essen.de\\
 {\bf Yuping Tu}\\[6pt]
Department of Mathematics, China Three Gorges University,\\
Yichang
443002, Hubei Province,  P. R. of China \\
tuyu02@sina.com\\[8pt]
 {\bf Alexis G. Zamora}\\[6pt]
CIMAT, AP 402, C.P. 36000, Guanajuato, Gto., M\'exico \\
alexis@cimat.mx }

%%%%%%%%%%% Introduction %%%%%%%%%%%%%%%%%%%%%%%%%%%%%%%%%%%
\def\d{\partial}
\def\c_t{\chi_{\text{top}}}
\def\BD{\boldsymbol D}
\def\BDelta{\boldsymbol\Delta}
\def\wt{\widetilde}
\def\wh{\widehat}
\def\la{\leftarrow}
\def\ep{\varepsilon}
\def\ra{\rightarrow}
\def\lra{\longrightarrow}
\def\G{\Gamma}
\def\O{\mathcal O}
\def\I{\mathcal I}
\def\L{\mathcal L}
\def\M{\mathcal M}
\def\E{\mathcal E}
\def\F{\mathcal F}
\def\H{\mathcal H}
\def\l{\ell}
\def\lra{\longrightarrow}
\def\ra{\rightarrow}
\def\Pic{\text{\rm{Pic}}}
\def\Spec{\text{\rm{Spec\,}}}
\def\grad{\text{\rm{grad\,}}}
\def\syz{\text{\rm{syz\,}}}
\def\tr{\text{\rm{tr\,}}}
\def\div{\text{\rm{div}}}
\def\Div{\text{\rm{Div}}}
\def\parlistitem#1{\par\hskip0.7cm\llap{(#1)}\ }
\def\ratmap{\,-\,-\to\,}
\def\phi{\varphi}
\def\({$($}
\def\){$)$}
\def\disc{\text{\rm disc\,}}

\font\eightrm=cmr10 at 8pt \font\bigrm=cmr10 scaled1440
\font\midrm=cmr10 scaled1200 \font\ninerm=cmr10 at 9pt
\font\ninebf=cmbx10 at 9pt \font\bbf=msbm10 \font\bigbf=cmbx10
scaled1440 \font\Bigbf=cmbx10 scaled2073 \font\midbf=cmbx10
scaled1200 \font\eightit=cmti8 \font\nineit=cmti9
%%%%%%%%%%%%%%%%%%%%%%%%%%%%%%%%%%%%%%%%%%%%%%%%%%%%%%%%%%%%%%%%%%%%%%%%%
%\newcommand{\dfrac}[2]{\frac{\displaystyle#1}{\displaystyle#2}}
\renewcommand{\labelenumi}{$($\arabic{enumi}$)$}
\renewcommand{\labelenumii}{(\Alph{enumi})}
\renewcommand{\atop}[2]{\genfrac{}{}{0pt}{}{#1}{#2}}
\date{}
\maketitle

\section*{Introduction}

We denote by $X$ a complex smooth projective surface, and by
$f:X\to C$ a fibration over a curve $C$ whose generic fiber is a
curve of genus $g$. $f$ is called {\it isotrivial} if all smooth
fibers are isomorphic to a fixed curve. $f$ is called {\it
semistable} if all of the singular fibers are reduced nodal
curves. If there is no $(-1)$-curve contained in the fibers, then
we call $f$ {\it relatively minimal}. The projection $X=F\times
C\to C$ is called a trivial fibration. $s$ is always the number of
singular fibers of $f$.

In a well-known paper \cite{Be1}, Beauville proved that a
non-isotrivial fibration $f:X\to \Bbb P^1$ admits at least 3
singular fibers, and if $f$ is semistable, then $s\geq 4$. The
first author proved in \cite{Ta1} that if $f$ is a semistable
fibration of genus $g\geq 2$, then it admits at least $5$ singular
fibers.

Beauville \cite{Be2} gave a beautiful explicit classification of
all semistable elliptic fibrations $(g=1)$ over $\Bbb P^1$ with
$4$ singular fibers. More precisely, there are exactly $6$
non-isotrivial such families, and all of them are modular families
of elliptic curves. A very interesting natural problem is to
classify semistable fibrations over $\Bbb P^1$ of genus $g\geq2$
with $5$ singular fibers. In the survey (\cite{Vie}, \S1 and
Ex.\,5.9), Viehweg considered this problem. The purpose of this
note is to try to get some information on the structure of the
surface $X$ when $f$ has $5$ or $6$ singular fibers. Our main
results are the following two theorems.

\begin{thm}\label{Theorem I} Assume that $f:X\to \Bbb P^1$ is a non trivial
semistable fibration of genus $g\geq 2$ with $s$ singular fibers.
Assume also that $f$ is relatively minimal. Then we have:
 \begin{enumerate}
 \item[\rm{(1)}] If $s=5$, then $X$ is birationally rational or ruled.
 \item[\rm{(2)}] If $s=6$ and $g=2,\ 3$ or $4$, then $X$ is not of general type.
 \item[\rm{(3)}] If $s=6$, $g=5$ and $X$ is of general type, then the
minimal model $S$ of $X$ satisfies
$$
K_S^2=1,\hskip0.3cm p_g(S)=2, \hskip0.3cm q(S)=0.
$$
The fibration $f$ comes from a pencil on $S$ with $5$ simple base
points $($including infinitely near base points$)$.
\end{enumerate}
\end{thm}

Thus if the Kodaira dimension of $X$ is non negative, then $s\geq
6$. We have no example of case (3), and we conjecture that $s$ is
at least 7 when $X$ is of general type. We shall give two examples
to show that the two bounds are sharp.

\begin{enumerate}
 \item[A)] $s=6$, $g=3$, $K_X^2=-4$, and $X$ is birationally
a K3 surface.

 \item[B)] $s=7$, $g=4$, $K_X^2=-3$, $p_g(X)=2$, $q(X)=0$, and $X$ is of
general type containing 4 $(-1)$-curves.
\end{enumerate}

 Our next purpose is to give some inequalities for genus
 $g\geq 2$ fibrations over $\Bbb P^1$ (not necessarily
 semistable), which imply the lower bounds on $s$.

 \begin{thm}\label{Theorem II} Let $f:X\to \Bbb P^1$ be a relatively minimal
 fibration of genus $g\geq 2$, and let $K_f=K_{X/\mathbb P^1}=K_X+2F$ be its
 relative canonical divisor. Assume that $f$ is not locally trivial. Then
  $$
  K_f^2\geq
  \begin{cases} 4g-4, & \text{ if } \ \kappa(X)=-\infty,\cr
  6g-6, & \text{ if } \ \kappa(X)=0, \cr
   6g-5, & \text{ if } \ \kappa(X)=1,\cr
6g-6+\dfrac12\left(K_S^2+\sqrt{K_S^2}\sqrt{K_S^2+8g-8}\right), &
\text{ if } \ \kappa(X)=2,\end{cases}
 $$
 where $S$ is the unique minimal model of $X$ when $\kappa(X)\geq 0$.
 \end{thm}
The bounds are optimal for infinitely many $g$. The example given
in \cite{Ta1} satisfies $K_f^2=4(g-1)$, $g=2$ and $s=5$. The
example A) above satisfies $K_f^2=6(g-1)$, $g=3$ and $s=6$, where
$X$ is a K3 surface. In Theorem \ref{Corollary 2}, we will
classify the fibrations $f:X\to \Bbb P^1$ with minimal $K_f^2$
according to its Kodaira dimension.

 The proof of Theorem \ref{Theorem I} is to use some
inequalities, particularly the following strict canonical class
inequality and its refinement for a non trivial semistable
fibration $f:X\to C$ of genus $g\geq 2$ with $s\neq 0$ \cite{Ta1}
(see also \cite{Liu} for a differential geometric proof):
\begin{equation}\label{0.1}
K_{X/C}^2<(2g-2)(2g(C)-2+s).%\tag0.1
\end{equation}
We use Reider's method to prove Theorem \ref{Theorem II}. Note
that the strict canonical class inequality and the inequality
$K_f^2\geq 6(g-1)$ (resp. $K_f^2\geq 4(g-1)$) in Theorem
\ref{Theorem II} imply that $s\geq 6$ (resp. $5$).
 \vskip0.4cm

\section{Preliminaries}

Let $f:X\to C$ be a fibration of genus $g\geq 2$, namely $X$
(resp. $C$) is a nonsingular complex surface (resp. curve) and the
generic fiber $F$ of $f$ is a nonsingular curve of genus $g$. We
always assume that $f$ is relatively minimal, i.e., there is no
$(-1)$-curve contained in the fibers. $f$ is called {\it
semistable} if all of the singular fibers are reduced nodal
curves. Denote by $K_f=K_{X/C}=K_X-f^*K_C$ the relative canonical
divisor of $f$, and by $\omega_{X/C}$ its corresponding invertible
sheaf. The relative invariants of $f$ are defined as follows:
 \begin{align*} & \chi_f=\deg
f_*\omega_{X/C}=\chi(\O_X)-(g-1)(g(C)-1), \cr &
K_f^2=K^2_{X/C}=K_X^2-8(g-1)(g(C)-1),\cr &
e_f=\chi_{\text{top}}(X)-4(g-1)(g(C)-1)=\sum_{i=1}^s(\chi_{\text{top}}(F_i)-(2-2g)).
\end{align*}
These invariants are nonnegative, and $K_f^2=0$ (equivalently,
$\chi_f=0$) if and only if $f$ is locally trivial. $e_f=0$ iff $f$
is smooth. Let
$$e_{F_i}=\chi_{\text{top}}(F_i)-(2-2g).$$
 If $F_i$ is semistable, then $e_{F_i}$ is equal to the number of nodes of $F_i$.

Arakelov \cite{Ara} and Beauville \cite{Be3} proved that $K_f$ is
a nef divisor and the curves $E$ with $EK_f=0$ are those
$(-2)$-curves in the fibers. The map defined by $|nK_f|$ for large
$n$ is the contraction morphism $\sigma:X\to X^\#$ of the vertical
$(-2)$-curves (cf. \cite{Ta2}), we get the relative canonical
model $f^\#: X^\#\to C$ of $f$.

In what follows, we assume that $f:X\to C$ is semistable with $s$
singular fibers. Denote by $q$ a singular point of $X^\#$. Then
$(X^\#,q)$ is a rational double point of type $A_{\mu_q}$, here
$\mu_q$ is the number of $(-2)$-curves in $X$ over $q$. We also
denote by $q$ the singular point of the fibers on the smooth part
of $X^\#$, in this case $\mu_q=0$. Then
$$e_f=\sum_q(\mu_q+1)$$
is the number of nodes in the singular fibers. For convenience, we
let
$$
r_f:=\sum_q\dfrac1{\mu_q+1}.
$$
It is obvious that
\begin{equation}\label{2.9}
  r_f\leq e_f.
\end{equation}
 In \cite{Ta1}, we proved the strict canonical
class inequality when $s\neq 0$:
$$
K_f^2<(2g-2)(2g(C)-2+s).
$$
The inequality follows from the following
\begin{equation}\label{2.1}
K_f^2-(2g-2)(2g(C)-2+s)\leq -\dfrac{(2g-2)s}{e}+\dfrac{3\cdot
r_f}{e^2},
\end{equation}
where $e\geq 2$ is any integer (see \cite{Ta1}, p.594, (5)).
Equivalently,
\begin{equation}\label{2.2}
  r_f\geq
  \dfrac13e^2\left(K_X^2-(2g-2)\left(6g(C)-6+s-\dfrac{s}{e}\right)\right).
\end{equation}

\begin{lemma}\label{Lemma 5} Let $q_1,\cdots,q_r$ be the points such that
$\mu_{q_i}\neq 0$. Let $\ell'=\sum_{\mu_q\neq 0}\mu_q$ be the
number of $(-2)$-curves contained in the fibers of $f$. Then
\begin{equation}\label{2.3}
r_f\leq e_f-\ell'-{r\over 2}.         %\tag2.3
\end{equation}
\end{lemma}
\begin{proof} Note that the number of points $q$ such that $\mu_q=0$
is $e_f-\sum_i(\mu_{q_i}+1)$. Hence
 \begin{align*} \sum_q{1\over \mu_q+1}
&=e_f-\sum_i(\mu_{q_i}+1)+\sum_i{1\over \mu_{q_i}+1}\cr &\leq
e_f-\ell'-r+{r\over 2}\cr &=e_f-\ell'-{r\over 2}.
\end{align*}
This completes the proof.
\end{proof}

\begin{lemma}\label{Lemma 6} $e_{F_i}=g-1+c_i-g(\wt F_i)$, where $c_i$ is
the number of components of $F_i$ and $\wt F_i$ is the
normalization of $F_i$. So
\begin{equation}\label{2.4}
e_f\leq s(g-1)+\sum_{i=1}^s c_i.                  %\tag2.4
\end{equation}
\end{lemma}
\begin{proof} Let $\sigma:\wt F_i\to F_i$ be the normalization of $F_i$, and
$\Delta$ be the subscheme of the singular locus of $F_i$. Then we
have
$$
0\to \O_{F_i}\to \sigma_*\O_{\wt F_i} \to \O_\Delta \to 0.
$$
Since $F_i$ has only nodes as its singular points, we have
$e_{F_i}=\deg\Delta=h^0(\O_\Delta)=\chi(\O_\Delta)$. Thus
 \begin{align*}
e_{F_i}&=\chi(\O_{\wt F_i})-\chi(\O_{F_i})\cr
       &=c_i-g(\wt F_i)-(1-g) \cr
       &\leq g-1+c_i.
\end{align*}
This completes the proof.
\end{proof}

\begin{corollary}\label{Corollary 7} Denote by $\ell=\sum_ic_i-\ell'$ the number
of curves in the singular fibers different from $(-2)$-curves.
Then we have
\begin{equation}\label{2.5}
r_f\leq s(g-1)+\ell-{r\over 2}.           %\tag2.5
\end{equation}
\end{corollary}
\begin{proof} This follows from (\ref{2.3}) and (\ref{2.4}).
\end{proof}

\begin{lemma}\label{Lemma 8} If $e^2>\frac{6g+3}{g-1}$, then
\begin{equation}\label{2.8}
K_f^2\leq \dfrac{2e(g-1)^2\left((2g(C)-2+s)e-s\right)}{e^2(g-1)-3(2g+1)}. %\tag2.8
\end{equation}
\end{lemma}
\begin{proof}
By Cornalba-Harris-Xiao's inequality \cite{CH, Xi},
$$
K_f^2\geq {4g-4\over g}\chi_f,
$$
we get
$$
e_f\leq {2g+1\over g-1}K_f^2.
$$
Then from (\ref{2.9}) and (\ref{2.1}), we obtain
$$
K_f^2-(2g-2)(2g(C)-2+s)\leq -{(2g-2)s\over e}+{3(2g+1)\over
e^2(g-1)}K_f^2,
$$
it implies (\ref{2.8}).
\end{proof}

\section{The proof of Theorem \ref{Theorem II}}

In this section, we let $f:X\to \Bbb P^1$ be a relatively minimal
fibration of genus $g\geq 2$. If $\kappa(X)\neq -\infty$, we
denote by $S$ the unique minimal model of $X$.

\begin{thm}\label{Corollary 2}  Assume that $f$ is not
     locally trivial. Then
$$
  K_f^2\geq
  \begin{cases} 4g-4, & \text{ if } \ \kappa(X)=-\infty,\cr
  6g-6, & \text{ if } \ \kappa(X)=0, \cr
   6g-5, & \text{ if } \ \kappa(X)=1,\cr
6g-6+\dfrac12\left(K_S^2+\sqrt{K_S^2}\sqrt{K_S^2+8g-8}\right), &
\text{ if } \ \kappa(X)=2.\end{cases}
 $$
 \begin{enumerate}
 \item[{\rm(1)}] $K_f^2=4g-4$ if and only if $X$ is the minimal resolution of
      the singularities of a double covering surface
      $Z\overset\pi\to \Bbb P^1\times C$ ramified
      over a curve of numerical type $2F_1+(2g+2-4g(C))F_2$,
      and the fibration $f$ is induced by the first projection $\mathbb P^1\times C\to\mathbb
      P^1$. Here $F_i$ is a fiber of the $i$-th projection of \, $\mathbb
      P^1\times C$.
 \item[{\rm(2)}] If $\kappa(X)=0$, then $K_f^2=6g-6$ iff  the fibration is
      induced by a pencil $\Lambda\subset|C|$ on its minimal model $S$ with $C^2=2g-2$
      simple base points $($including infinitely near base points$)$.
 \item[{\rm(3)}] If $\kappa(X)=1$, then
      $K_f^2= 6g-5$ iff the fibration is
      induced by a pencil $\Lambda\subset|C|$ of genus $g$ on its minimal model $S$ with $C^2=2g-3$
      simple base points.
 \item[{\rm(4)}] If $\kappa(X)=2$, then $K_S^2\geq 1$.
      $K_f^2=
      6(g-1)+\dfrac12\left(K_S^2+\sqrt{K_S^2}\sqrt{K_S^2+8g-8}\right)$
      iff \ the fibration is induced by a pencil $\Lambda\subset |C|$ on $S$ with only simple
      base points and $C\sim rK_S$.
          In particular, if $g=2$, then $K_f^2\geq 8$, with equality iff $K_S^2=1$ and $C\sim K_S$.
      If $g=3$, then $K_f^2\geq 15$.
\end{enumerate}
\end{thm}
\begin{proof}
(1) Let $A=K_X+F$. We claim that $|2A|$ (resp. $|3A|$) is base
    point free if $g\geq 3$ (resp. $g=2$). In particular, $A$ is nef
    and hence
  $$
  K^2_f-4(g-1)=A^2\geq 0.
  $$

    Indeed, by Ramanujan's vanishing theorem (\cite{BPV}, p.131),
    $H^1(-F)=H^0(-F)=0$. So
  $$
  h^0(A)=h^2(-F)=\chi(-F)=\chi_f>0.
  $$
    Hence we can assume that $A$ is an effective divisor.  Since
    $AF=2g-2\geq 2$, $A$ admits at least one horizontal component.
    Hence $AK_f\geq 1$.

     Note that $|2A|=|K_X+K_f|$. If $g\geq 3$, $K_f^2=AK_f+2(g-1)\geq
     1+4=5$. Suppose $|2A|$ has a base point $p$, by Reider's theorem
     \cite{Rei}, there is a curve $E$ passing through $p$ such that

       (i) \ $K_fE=0$ and $E^2=-1$; \ \ or \ \ (ii) \ $K_fE=1$ and $E^2=0$.
                      \newline
    On the other hand, $K_fE=K_XE+2FE\equiv K_XE\equiv E^2 \pmod 2$.
    Thus the two cases can not exist. This proves that $|2A|$ has no
    base point.

    Now we consider the case $g=2$. We first prove that $A$ is nef.

    Suppose $\Gamma$ is an irreducible and reduced curve with $A\Gamma<0$.
    It is easy to see that $\Gamma$ must be a horizontal curve. Hence
    $(A+F)\Gamma=K_f\Gamma\geq 1$. Then we have $\Gamma F\geq 2$. On
    the other hand, $|A|$ is non empty, so $\Gamma$ is the fixed part
    of $|A|$. Hence $2=2g-2=AF\geq \Gamma F\geq 2$. We get $\Gamma
    F=2$, and $A=\Gamma+E$, where $E$ consists of vertical curves. It
    implies that $\Gamma^2=A\Gamma -E\Gamma<0$. Since
    $A\Gamma=K_X\Gamma+2<0$, we have $K_X\Gamma\leq -3$. So
    $p_a(\Gamma)<0$, impossible. Therefore, $A$ is nef and $A^2\geq
    0$.

    In this case, $|3A|=|K_X+2A+F|$. $L=2A+F$ is nef and
    $L^2=4A^2+4AF\geq 8(g-1)>4$. If $|3A|$ has a base point $p$, then
    there is a curve $E$ passing through $p$ such that

        (i) \ $LE=0$, $E^2=-1$; \ \ or \ \ (ii) $LE=1$, $E^2=0$.
                              \newline
    Since $A$ and $F$ are nef, (i) implies $AE=FE=0$. Note that
    $A=K_X+F$, we have $K_XE=0$ and $K_XE\not\equiv E^2$, a
    contradiction. So case (i) is impossible. In case (ii), we have
    $AE=0$ and $FE=1$, so $K_XE=-1$ which is also impossible since
    $E^2=0$. This proves that $|3A|$ is base point free.

    Now we consider the case when $K^2_f=4g-4$, i.e., $A^2=0$.

    In this case, the base point free linear system $|6A|$ is composed
    with a fibration $\varphi:X\to C$ over a smooth curve $C$ . Denote
    by $F'$ a generic fiber of $\varphi$. Then $AF'=0$. Since
    $AF=2g-2\geq 2$, $F'$ can not be the fiber of $f$. So $FF'\geq 1$.
    From $AF'=0$, we get $K_XF'=-FF'\leq -1$. Since ${F'}^2=0$, we
    have $K_XF'=-2$ and $FF'=2$. So $F'\cong \mathbb P^1$. Now we
    get a generally double cover $\pi:X\to \mathbb P^1\times C$
    defined by $\pi(x)=(f(x),\varphi(x))$. Since the pullback
    $\pi^*F_2$ of a generic fiber $F_2$ is still isomorphic to $\mathbb
    P^1$, and a generic $\pi^*F_1$ has genus $g$. By Hurwitz formula, we
    see that the branch curve $B_\pi$ of $\pi$ has numerical type
    $2F_1+(2g+2-4g(C))F_2$.

    Conversely, if the branch curve $B_\pi$ has numerical type
    $2F_1+(2g+2-4g(C))F_2$, then $B_\pi$ is nonsingular or admits at
    most $ADE$ singularities. Thus the canonical resolution is the
    minimal one, and the fibration induced is relatively minimal.
    By easy computations (\cite{BPV}, p.183, or \cite{Hor}), we have $K_X^2=-4(g-1)$
    and hence $K_f^2=4g-4$. This completes the proof of (1).

(2) Suppose that the Kodaira dimension of $X$ is non negative, and
    let $S$ be the unique minimal model of $X$ obtained by contracting
    $(-1)$-curves. So the fibers $F$ of $f$ are contracted to a pencil
    $\Lambda$ in $|C|$ on $S$ with base points $p_1, \cdots, p_m$
    $($including infinitely near base points$)$. We consider a generic
    curve $C$ in $\Lambda$, and let $n_i$ be the multiplicity of $C$
    at $p_i$. Then $n_i\geq 1$. Note that $F^2=0$, $K_XF=2g-2$. One
    can prove easily that
    \begin{alignat}{2}
    C^2&=F^2+\sum_{i=1}^mn_i^2= \sum_{i=1}^mn_i^2.&\label{1.1}
    \\
    K_SC&=K_XF-\sum_{i=1}^mn_i =2g-2-\sum_{i=1}^mn_i.&\label{1.2}
    \end{alignat}
    Since $S$ is neither ruled nor rational, $K_S$ is nef and
    $K_S^2\geq 0$. Hence $K_SC\geq 0$. Thus $m\leq \sum_in_i\leq
    2g-2-K_SC$.
     \begin{equation}\label{1.3}
     K_f^2=K_X^2+8(g-1)=K_S^2-m+8(g-1)\geq K_S^2+K_SC+6(g-1). %\tag1.3
     \end{equation}
    Hence $K_f^2\geq 6g-6$. If $K_f^2=6g-6$, i.e., $K_X^2=-(2g-2)<0$,
    then $K_SC=0$, $K_S^2=0$ and $n_i=1$ for any $i$. In particular,
    $m=K_S^2-K_X^2=2g-2>0$ and $C^2>0$. By Hodge index theorem, $K_S$
    is numerically trivial, thus $\kappa(X)=0$. By (\ref{1.1}),
    $C^2=m=2g-2$. This completes the proof of (2).

 (3) We have proved that if $\kappa(X)\geq 0$ and $K_f^2=6g-6$, then
    $\kappa(X)=0$. So if $\kappa(X)=1$, we obtain $K_f^2\geq 6g-5$ and $K_S^2=0$.
    Suppose $K_f^2=6g-5$, i.e., $K_X^2=-(2g-3)<0$, then $m=2g-3>0$ and $C^2>0$.
    As in the proof of (2), $K_SC\neq 0$ (otherwise, by Hodge
    index theorem, $K_S\sim 0$ and $\kappa(X)=0$). So $K_SC\geq
    1$.  From our assumption $K_f^2=6g-5$ and (2.3), we get $K_SC=1$.  Now
    (\ref{1.2}) implies that $n_i=1$
    for any $i$, hence $C^2=2g-3$ and $C$ is smooth. The fibration
    is induced by a pencil in $|C|$ with $2g-3$ simple base
    points.

(4) In the case when $\kappa(X)=2$, we let $x=CK_S$, $y=C^2$. Then
    \begin{equation}\label{1.4}
      x+y\geq 2g-2,\hskip0.3cm K_S^2\cdot y\leq x^2. %\tag1.4
    \end{equation}
    It is easy to prove that
  $$
  x\geq \dfrac12\left(-K_S^2+\sqrt{K_S^2}\sqrt{K_S^2+8g-8}\right).
  $$
     Now from (\ref{1.3}), we have
  $$
    K_f^2\geq
    6g-6+\dfrac12\left(K_S^2+\sqrt{K_S^2}\sqrt{K_S^2+8g-8}\right).
  $$
    If the equality holds, then
    the equalities in (\ref{1.4}) hold. We obtain easily the desired
    characterization.
\end{proof}

\begin{corollary}\label{Theorem 1} Let $f:X\to \Bbb P^1$ be a relatively
minimal semistable fibration of genus $g\geq 2$ with $s$ singular
fibers. If $f$ is non trivial and $s=5$,  then $X$ is birationally
ruled or rational.
\end{corollary}
\begin{proof} By the strict canonical class inequality, we have
$K_f^2<6g-6$. The corollary follows from Theorem \ref{Corollary
2}.
\end{proof}

Note that if $s=5$, then $S$ is either $\Bbb P^2$ or a
geometrically ruled surface. So
$$
m=K_S^2-K_X^2>K_S^2+2g-2=\begin{cases} 2g+7, &\text{ if }
S\cong\Bbb P^2,\cr 2g+6-8q(X), &\text{ otherwise.}
\end{cases}
$$

\section{The proof of Theorem \ref{Theorem I} for $s=6$}

In this section, we assume that $f:X\to\Bbb P^1$ is semistable and
$s=6$. We use freely the notations in the previous sections,
including those in the proof. Denote by $F_1,\cdots,F_6$ the 6
singular fibers.

We recall the inequalities in Sect.~1 for our $f$:
 \begin{alignat*}2
   &{\rm (A)} \hskip0.5cm K_X^2<0, \hskip8cm &\textrm{by (\ref{0.1})}\cr
   &{\rm (B)} \hskip0.5cm r_f \geq \dfrac13 e\left(K_X^2e+12g-12\right),  &\textrm{by (\ref{2.2})}\cr
   &{\rm (C)} \hskip0.5cm r_f\leq 6g-6+\ell-\dfrac{r}2, &\textrm{by
   (\ref{2.5})}\cr
   &{\rm (D)} \hskip0.5cm K_X^2\leq
   4\dfrac{(g-1)^2(2e-3)e}{e^2(g-1)-6g-3}-8(g-1), \hskip0.3cm \textrm{ if }
  \  e^2>\dfrac{6g+3}{g-1},      &\textrm{by (\ref{2.8})}
 \end{alignat*}
where $r_f=\sum_q{1\over \mu_q+1}$, and $\ell$ is the number of
non $(-2)$-curves in the 6 singular fibers.

In what follows, we assume that $X$ is of general type, thus
$K_S^2\geq 1$. Denote respectively by $C_1,\cdots,C_6\subset S$
the image curves of $F_1,\cdots,F_6$. By (A),
\begin{equation}\label{2.6}
m=K_S^2-K_X^2\geq 2.                 %\tag2.6
\end{equation}
We see that $K_SC\neq 1$. Otherwise, by Hodge index theorem, $
K_S^2C^2\leq (K_SC)^2=1, $ so $C^2\leq 1$, which implies $m\leq
1$. Hence
\begin{equation}\label{2.7}
2\leq K_SC\leq 2g-2-m.                   %\tag2.7
\end{equation}

\begin{thm}\label{Theorem 3} If $g\leq 4$, then $X$ is not of
general type.
\end{thm}
\begin{proof} We shall rule out the cases $g=2,3$ and $4$
under the assumption $\kappa(X)=2$.

 {\it Case} \ $g=2$: (\ref{2.7}) implies $m=0$, which contradicts (\ref{2.6}).

{\it Case} \ $g=3$: \ By (\ref{1.2}), (\ref{2.6}) and (\ref{2.7}),
$m=K_SC= 2$, $K_S^2=1$, $K_X^2=-1$, $n_1=n_2=1$, hence $C^2=2$. By
Noether inequality $K_S^2\geq 2p_g(S)-4$ and Debarre's inequality
\cite{Deb} $K_S^2\geq 2p_g(S)\geq 2q(S)$ if $q(S)\geq 1$, we
obtain
\begin{equation}\label{2.10}
p_g(S)\leq 2, \hskip0.3cm q(S)=0.     %\tag2.10
\end{equation}
Note that
\begin{equation}\label{2.11}
e_f=12\chi_f-K_f^2
   =12(p_g(S)+1)-K_X^2+4(g-1)
   =12p_g(S)+21\leq 45. %\tag2.11
\end{equation}
On the other hand, from (B) for $e=12$, we get $ r_f\geq 48.$ This
contradicts (\ref{2.9}): $r_f\leq e_f$.

{\it Case} \ $g=4$:  \ (D) for $e=12$ implies that $K_X^2\leq
-8/5$. Thus $K_X^2\leq -2$ and $m\geq 3$. By (\ref{2.7}), we see
that $m\leq 4$.

If $m=3$, then $K_S^2=1$ and $K_X^2=-2$, we see that (\ref{2.10})
holds.  As in (\ref{2.11}), we obtain $e_f\leq 50$. (\ref{2.2})
for $e=9$ implies that $r_f\geq 54$, which contradicts $r_f\leq
e_f$.

If $m=4$, then (\ref{1.2}) and (\ref{2.7}) imply that $K_SC=2$,
$n_1=\cdots=n_4=1$, and $C^2= 4$. By Hodge index theorem, we can
see that $K_S^2=1$ and $K_S^2C^2=(K_SC)^2$. Hence $C\sim 2K_S$.
This means that any $(-2)$-curve $E$ on $S$ does not pass through
any one of the 4 base points $p_1,\cdots,p_4$ because
$CE=2K_SE=0$. Hence $(-2)$-curves in $C_i$ must be $(-2)$-curves
in $F_i$.

 Note that $K_SC_i=2$, so $C_i$ has at most two
components different from $(-2)$-curves. This implies that $F_i$
has at most two components which are neither $(-2)$-curves nor the
exceptional curves of the 4 base points.

On the other hand, among the 4 exceptional curves on $X$, at least
one of them is a horizontal $(-1)$-curve. If the remaining three
exceptional curves are vertical, then at least one of them is a
$(-2)$-curve. Hence the total number $\ell$ of non $(-2)$-curves
in the 6 singular fibers $\leq 6\times 2+2=14$.

Now by (C), $r_f\leq 32$. On the other hand, (B) for $e=6$ implies
that $r_f\geq 36$, a contradiction.
\end{proof}

\begin{thm}\label{Theorem 9} If $s=6$, $g=5$ and $X$ is of general type,
then the minimal model $S$ of $X$ satisfies
$$
K_S^2=1,\hskip0.3cm p_g(S)=2, \hskip0.3cm q(S)=0.
$$
The fibration comes from a pencil on $S$ with $5$ simple base
points.
\end{thm}
\begin{proof} In this case, (D) for $e=10$ implies that $
K_X^2\leq -{864\over 367}<-2, $ thus $K_X^2\leq -3$ and $m\geq 4$.
On the other hand, $m\leq 2g-2-2=6$. If $m=6$, then $C^2\geq 6$
and $K_SC\leq 2$, by Hodge index theorem,
$K_S^2C^2\leq(K_SC)^2\leq 4$, and hence $K_S^2\leq 0$, a
contradiction. Thus $m=4$ or $5$.

If $m=4$, then $K_S^2=1$, $K_X^2=-3$, $p_g(S)\leq 2$, $q(S)=0$,
 $C^2\geq 4$, and $K_SC\leq 4$. Thus
$$
e_f=12(p_g+1)-K_X^2+4(g-1)\leq 55.
$$
From (\ref{2.2}) for $e=8$, we have
$$
r_f\geq {1\over 3}e(-3e+48)=e(16-e)=64>e_f.
$$
It contradicts (\ref{2.9}).

So $m=5$. By (\ref{1.1}) and (\ref{1.2}), we have $C^2\geq 5$,
$K_SC\leq 3$. From (\ref{2.7}), we get $K_SC\geq 2$. Since
$K_S^2\geq 1$, by Hodge index theorem, $K_S^2C^2\leq (K_SC)^2$, we
have $K_SC=3$ and $K_S^2=1$. So $K_X^2=-4$, $p_g\leq 2$ and $q=0$.
By (\ref{1.1}) and (\ref{1.2}), we obtain $n_1=\cdots=n_5=1$ and
$C^2=5$. Generic $C$ in the pencil is smooth. Let $e=6$. Then (B)
gives us that $r_f\geq {1\over3}e(-4e+48)=48$.

Now we prove that $p_g=2$. Otherwise $p_g\leq 1$, by (\ref{2.11}),
$e_f\leq 12(1+1)-(-4)+16=44<r_f$, a contradiction. Then we have
$$
K_S^2=1,\hskip0.3cm p_g=2, \hskip0.3cm q=0, \hskip0.3cm m=5.
$$
So $e_f=\sum_q(\mu_q+1)=56$, $r_f=\sum_q\frac1{\mu_q+1}\geq 48$.

Let $m_i=\#\{\,q\,|\,\mu_q=i\,\}.$ Then
 \begin{align*}
e_f&=m_0+2m_1+3m_2+\cdots=56, \cr
 r_f&=m_0+\dfrac12 m_1+\dfrac13 m_2+\cdots\geq 48.
\end{align*}
We see that $e_f-m_0\geq 4(r_f-m_0)$, hence $m_0\geq 46$. If
$m_0=46$, one can prove that $m_1=5$ and $m_i=0$ for $i>1$. So
there are at least 46 points $q$ such that $\mu_q=0$, and there
are at most $(56-46)/2=5$ points $q$ with $\mu_q\neq 0$. By
(\ref{2.3}) we get
$$
\ell'+{r\over 2}\leq 8.
$$
It implies that $r\leq 5$ and $\ell'\leq 7$.
\end{proof}

\section{Two examples}

Note that for any two points $p_1,p_2$ on $\Bbb P^1=\Bbb
C\cup\{\infty\}$, there exists a cyclic cover $\phi_{p_1,p_2}:\Bbb
P^1\to \Bbb P^1$ ramified exactly over $p_1$ and $p_2$. For
example, $\tau_n:\Bbb P^1\to \Bbb P^1$ defined by $x\mapsto x^n$
is totally ramified over $0$ and $\infty$. Now we are going to
construct some covers of $\Bbb P^1$.

Let $\pi:\Bbb P^1\to\Bbb P^1$ be the double cover ramified over $1$ and
$\infty$, and let $\psi_e:\Bbb P^1\to \Bbb P^1$ be the cyclic cover of
degree $e$ totally ramified over the two points $\pi^{-1}(0)$.
Then the composition
$$
\phi_{2e}:\Bbb P^1 \ \overset{\psi_e}\to \ \Bbb P^1 \
\overset{\pi}\to\ \Bbb P^1
$$
is a covering of degree $2e$ ramified uniformly over $0,1,\infty$.
The map $\varphi_{2e}$ is the quotient map of $\mathbb P^1$ for
the standard action of the diedral group of order $2e$ on $\mathbb
P^1$.
 The cover has two ramification points over $0$
with ramification index $e$, and $e$ ramification points over $1$
(resp. $\infty$) with index 2. Now we consider the fiber product
of $\tau_n:\Bbb P^1\to \Bbb P^1$ and $\phi_{2e}: \Bbb P^1\to\Bbb
P^1$, we get a curve
$$
\Gamma=\{\,(x,y)\in \Bbb P^1\times\Bbb P^1 \,|\, \tau_n(x)=\phi_{2e}(y)\,\}.
$$
This curve is of type $(2e,n)$. The first projection $p_1:\Gamma
\to \Bbb P^1$ can be viewed as the pullback of $\phi_{2e}$ under
the base change $\tau_n$. So the set $\Sigma$ of critical points
of $p_1$ is $\tau_n^{-1}\{0,1,\infty\}$. $\Sigma$ contains $n+2$
points. In fact, $\Gamma $ is locally defined by
$\phi_{2e}(y)=\tau_{n}(x)$. Thus we can see that $\Gamma$ admits
two singular points defined by $y^e=x^n$ over $\tau_n^{-1}(0)$ and
$e$ singular points of type $y^2=x^n$ over $\tau_n^{-1}(\infty)$.
$p_1$ has $e$ simple ramification points over each point of
$\tau_n^{-1}(1)$.
$$
\begin{CD} \Gamma\hskip0.3cm \subset \hskip0.3cm @.\mathbb P^1\times\mathbb
P^1@>\textrm{pr}_2>>\mathbb P^1\\
p_1\searrow@.@V{\textrm{pr}_1}VV@V{2e:1}V{\phi_{2e}}V \\
@.\mathbb P^1@>{n:1}>{\tau_n}>\mathbb P^1
\end{CD}
$$

Let $e=n=4$. Then $\Gamma$ is of type $(8,4)$. Let $X$ be the
minimal resolution of the double cover $\Sigma\to \Bbb
P^1\times\Bbb P^1$ ramified over $\Gamma$. $X$ can be obtained by
the canonical resolution $\sigma:Y\to \mathbb P^1\times\Bbb P^1$,
where $\sigma$ is the blowing-up of $\mathbb P^1\times\Bbb P^1$ at
the two singular points of $\Gamma$ over $\tau_n^{-1}(0)$. Denote
by $E_1$ and $E_2$ the two exceptional curves, and by $\bar\Gamma$
the strict transform of $\Gamma$ on $Y$. Then $\bar\Gamma \equiv
\sigma((8,4))-4E_1-4E_2=2\bar\delta$, where
$\bar\delta=\sigma^*((4,2))-2E_1-2E_2$. Let $\bar\pi:\bar X\to Y$
be the double cover ramified over $\bar\Gamma$. Note that
$\bar\Gamma$ admits only ADE singularities, so $\bar X$ admits at
most rational double points, and $X$ is the minimal resolution of
$\bar X$. Now by the formulas (\cite{BPV}, p.183), we have
$$
K_{\bar
X}=\bar\pi^*(K_Y+\bar\delta)=\bar\pi^*(\sigma^*(2,0)-E_1-E_2).
$$
Denote by $C_1$ and $C_2$ the two horizontal sections of
$\text{pr}_1:\mathbb P^1\times\mathbb P^1 \to \mathbb P^1$ passing
through the two points blown up by $\sigma$, and by $\bar C_i$
their strict transforms on $Y$. Since $\bar C_1$ and $\bar C_2$
are disjoint with $\bar \Gamma$, it is easy to see that their
pullback on $\bar X$ are 4 $(-1)$-curves. Note that
$$\bar C_1+\bar C_2=\sigma^*(2,0)-E_1-E_2,$$
so $K_{\bar X}=\bar\pi^*(\bar C_1+\bar C_2)$ is the sum of the 4
$(-1)$-curves. Let $S$ be the surfaces obtained by contracting the
4 horizontal $(-1)$-curves. Then we see that $K_S\equiv 0$. By
easy computation, we have
\begin{align*}
\chi(\mathcal O_{S})&=\chi(\mathcal O_{\bar X})\\
                         &=2\chi(\mathcal
                         O_Y)+\dfrac12(\bar\delta^2+K_Y\bar\delta)\\
                         &=2+0=2.
\end{align*}
 So $S$ is a $K3$ surface. The first
projection of $\Bbb P^1\times\Bbb P^1$ induces a semistable
fibration $f:X\to \Bbb P^1$ of genus $3$ with $n+2=6$ singular
fibers. Because $K_X^2=-4$, $K_f^2=12=6g-6$.

Similarly, let $e=n=5$, and let $F_\infty$ be the vertical fiber
of $\Bbb P^1\times \Bbb P^1$ passing through the two singular
points of type $x^5=y^5$. Then consider the double cover of $\Bbb
P^1\times \Bbb P^1$ ramified over $F_\infty +\Gamma$. The
canonical resolution is similar to the above case. We have $\bar
\delta \equiv \sigma^*((5,3))-3E_1-3E_2$ and
\begin{align*}
K_Y+\bar\delta &=\sigma^*((3,1))-2E_1-2E_2\\
               &=\sigma^*((1,1))-E_1-E_1+\bar C_1+\bar C_2\\
               &=\sigma^*((1,0))+\bar F_\infty+\bar C_1+\bar C_2,
\end{align*}
where $\bar F_\infty\equiv\sigma^*((0,1))-E_1-E_2$ is the strict
transform of $F_\infty$ on $Y$. The pullback of $\bar
F_\infty+\bar C_1+\bar C_2$ on $\bar X$ consists 5 $(-1)$-curves,
the one from $\bar F_\infty$ is a vertical $(-1)$-curve. After
blowing-down this vertical $(-1)$-curve, we get a relatively
minimal semistable fibration $f:X\to \mathbb P^1$ of genus $g=4$
with $n+2=7$ singular fibers. We see easily that $|K_S|$ is a
genus $2$ pencil with one simple base point. The pencil is in fact
induced by pr$_2$. So $X$ is of general type. By easy computation,
we have $\chi(\mathcal O_S)=3$, hence
  $$
    K_S^2=1, \ p_g(S)=2, \ q(S)=0.
  $$
In this case, $K_X^2=-3$ and $K_f^2=6g-3$.

\vskip0.3cm
 \noindent{\bf Acknowlegment.} The first author would like to thank
 Prof. H. Esnault and Prof. E. Viehweg for their hospitality and
 useful discussions with Viehweg while he is visiting Universt\"at Essen.
 The third author would like to thank Prof. A. Beauville for
 useful comments. We all thank the referee for pointing out a
 mistake and providing us many valuable suggestions in rewriting the
 manuscript which make the paper readable.


\begin{thebibliography}{99}

\bibitem{Ara} S. Ju. Arakelov: {\it Families of algebraic curves with
fixed degeneracy}, Math. USSR Izvestija {\bf 5} (1971), 1277--1302

\bibitem{BPV}
 W. Barth, C. Peters, A. Van de Ven: Compact Complex Surfaces, Springer Verlag
   1984


\bibitem{Be1} A. Beauville: {\it Le nombre minimum de fibres
singuli\`eres d'un courbe stable sur $\Bbb P^1$,} in: S\'eminaire
sur les pinceaux de courbes de genre au moins deux (L. Szpiro,
ed.), Ast\'erisque  {\bf 86} (1981), 97--108



\bibitem{Be2} A. Beauville: {\it Les familles stables de courbes
elliptiques sur $\Bbb P^1$ admettant quatre fibres singuli\`eres},
 C.R. Acad. Sci. Paris {\bf 294} (1982), 657--660

\bibitem{Be3} A. Beauville: {\it L'in\'eqalit\'e $p_g\geq 2q-4$ pour
 les surfaces de type g\'en\'eral}, Appendice \'a O. Debarre:
 {\it In\'eqalit\'es num\'eriques pour les surfaces de type
 g\'en\'eral}, Bull. Soc. Math. France {\bf 110} (1982), no.~3, 319--346

\bibitem{CH}
 M. Cornalba, J. Harris:
 {\it Divisor
classes associated to families of stable varieties, with
applications to the moduli space of curves}, Ann. Sci. \'Ecole
Norm. Sup. (4) {\bf 21} (1988) no.~3, 455--475

\bibitem{Deb}
O. Debarre: {\it In\'eqalit\'es num\'eriques pour les surfaces de
type g\'en\'eral}, Bull. Soc. Math. France {\bf 110} (1982),
 no.~3, 319--346


\bibitem{Hor} E. Horikawa:
 {\it Algebraic surfaces of general type
 with small $c^{2}_{1}$. V}
   J. Fac. Sci. Univ. Tokyo Sect. IA Math.
  {\bf 28} (1981), no.~3, 745--755



\bibitem{Liu}
 K. Liu:
 {\it Geometric height inequalities},
 Math. Res. Lett. {\bf 3} (1996) no.~5, 693--702



%\bibitem{Kon} K. Konno:
% {\it Nonhyperelliptic fibrations of
%small genus and certain irregular canonical surfaces},
%  Ann. Scuola Norm. Sup. Pisa Cl. Sci.
%   {\bf 20} (1993) no.~4, 575--595


\bibitem{Rei}
 I. Reider:
 {\it Vector bundles of rank $2$ and linear systems on algebraic
 surfaces}
  Ann. of Math.
  {\bf 127} (1988), 309--316


\bibitem{Ta1} S.-L. Tan: {\it The minimal number of singular
fibers of a semistable curve over $\Bbb P^1$},  J. Algebraic
Geometry   {\bf 4} (1995), 591--596


\bibitem{Ta2}
 S.-L. Tan:
 {\it Effective behavior of multiple linear systems},
  Asian J. of Math.
  {\bf 7} (2003) no.~4


\bibitem{Vie} E. Viehweg:
 {\it Positivity of direct image
sheaves and applications to families of higher dimensional
manifolds}, in:
 School on Vanishing Theorems and Effective
Results in Algebraic Geometry (Trieste, 2000)  249--284
  ICTP Lect. Notes, 6, Int. Cent. Theoret. Phys., Trieste,
2001


\bibitem{Xi}
 G. Xiao:
 {\it Fibered algebraic surfaces with low slope},
  Math. Ann. {\bf 276} (1987)
 no.~3, 449--466


\end{thebibliography}
\enddocument